\documentclass[10pt]{article}
\usepackage{amsfonts, amsmath,amssymb}

\newcommand{\detail}[1]{\par\noi{\bf [Proof detail\ }{#1}
\hfill{\bf ]}\par\noi\hspace{-4pt}}
\renewcommand{\detail}[1]{}

\newcommand{\file}{papers/braco/errata.tex\quad}
\renewcommand{\file}{}

\newcommand{\dis}{\displaystyle}

\newcommand{\noi}{\noindent}

\newcommand{\halmos}{\rule{1ex}{1.4ex}}
\def \qed {\nopagebreak{\hspace*{\fill}$\halmos$\medskip}}

\newtheorem{theorem}{Theorem}
\newtheorem{proposition}[theorem]{Proposition}
\newtheorem{corollary}[theorem]{Corollary}
\newtheorem{conjecture}[theorem]{Conjecture}
\newtheorem{lemma}[theorem]{Lemma}

\newtheorem{counterexample}[theorem]{Counterexample}
\newtheorem{remark}[theorem]{Remark}

\newcommand{\bt}{\begin{theorem}}
\newcommand{\et}{\end{theorem}}
\newcommand{\bl}{\begin{lemma}}
\newcommand{\el}{\end{lemma}}
\newcommand{\bp}{\begin{proposition}}
\newcommand{\ep}{\end{proposition}}
\newcommand{\bcor}{\begin{corollary}}
\newcommand{\ecor}{\end{corollary}}
\newcommand{\br}{\begin{remark}\rm}
\newcommand{\er}{\end{remark}}
\newcommand{\bcon}{\begin{conjecture}}
\newcommand{\econ}{\end{conjecture}}

\newcommand{\be}{\begin{equation}}
\newcommand{\ee}{\end{equation}}

\newcommand{\bes}{\begin{equation*}}
\newcommand{\ees}{\end{equation*}}
\newcommand{\bea}{\begin{eqnarray}}
\newcommand{\eea}{\end{eqnarray}}
\newcommand{\ba}{\begin{array}}
\newcommand{\ea}{\end{array}}
\newcommand{\bc}{\be\begin{array}{r@{\,}c@{\,}l}}
\newcommand{\ec}{\end{array}\ee}

\newcommand{\oo}{\omega}
\newcommand{\om}{\Omega}

\newcommand{\Di}{{\cal D}}

\newcommand{\Fi}{{\cal F}}

\newcommand{\R}{{\mathbb R}}

\newcommand{\sub}{\subset}

\newcommand{\asto}[1]{\underset{{#1}\to\infty}{\longrightarrow}}

\newcommand{\ffrac}[2]{{\textstyle\frac{{#1}}{{#2}}}}
\newcommand{\dif}[1]{\ffrac{\partial}{\partial{#1}}}

\newcommand{\di}{\mathrm{d}}
\newcommand{\half}{{[0,\infty)}}

\newcommand{\spa}{{\rm span}}

\begin{document}

\makeatletter\@addtoreset{equation}{section}
\makeatother\def\theequation{\thesection.\arabic{equation}} 

\renewcommand{\labelenumi}{{(\roman{enumi})}}

\title{Correction to: Branching-coalescing particle systems}
\author{Siva R.~Athreya \and Jan M. Swart}
\date{\file\today}
\maketitle

\begin{abstract}\noi
In the article titled ``Branching-Coalescing Particle Systems'' published in
{\em Probability Theory and Related Fields~131(3), pages 376--414, (2005)},
Theorem~7 as stated there is incorrect. Indeed, we show by counterexample that
the equality that we claimed there to hold for all time, in general holds only
for almost every time with respect to Lebesgue measure. We prove a weaker
version of the theorem that is still sufficient for our applications in the
mentioned paper.
\end{abstract}
\vspace{20pt}

\noi
{\it Acknowledgements.} Siva Athreya is supported in part by a CSIR grant. Jan
Swart is supported by GA\v CR grant 201/09/1931. We thank Martin Hutzenthaler
for pointing out the error in \cite[Theorem~7]{AS05}.
\vspace{12pt}

\section{Introduction}

In this note we repair an error in \cite[Theorem~7]{AS05}. We start by
stating the corrected theorem.

If $E$ be a metrizable space, we denote by $M(E),B(E)$ the spaces of
real Borel measurable and bounded real Borel measurable functions on
$E$, respectively. If $A$ is a linear operator from a domain
$\Di(A)\sub M(E)$ into $M(E)$ and $X$ is an $E$-valued process, then
we say that $X$ solves the martingale problem for $A$ if $X$ has
cadlag sample paths and for each $f\in\Di(A)$,
\be\label{MP}
E\big[|f(X_t)|\big]<\infty\quad\mbox{and}\quad
\int_0^tE\big[|Af(X_s)|\big]\di s<\infty\qquad(t\geq 0),
\ee
and the process $(M_t)_{t\geq 0}$ defined by
\be\label{duM}
M_t:=f(X_t)-\int_0^t\!Af(X_s)\di s\qquad(t\geq 0)
\ee
is a martingale with respect to the filtration generated by $X$.

We will prove the following theorem.
\bt{\bf(Duality with error term)}\label{errorth}
Assume that $E_1,E_2$ are metrizable spaces and that for $i=1,2$,
$A_i$ is a linear operator from a domain $\Di(A_i)\sub B(E_i)$ into
$M(E_i)$. Assume that $\Psi\in B(E_1\times E_2)$ satisfies
$\Psi(\cdot,x_2)\in\Di(A_1)$ and $\Psi(x_1,\cdot)\in\Di(A_2)$ for each
$x_1\in E_1$ and $x_2\in E_2$, and that
\be
\Phi_1(x_1,x_2):=A_1\Psi(\,\cdot\,,x_2)(x_1)\quad\mbox{and}\quad\Phi_2(x_1,x_2)
:=A_2\Psi(x_1,\,\cdot\,)(x_2)\qquad(x_1\in E_1,\ x_2\in E_2)
\ee
are jointly measurable in $x_1$ and $x_2$. Assume that $X^1$ and $X^2$
are independent solutions to the martingale problems for $A_1$ and
$A_2$, respectively, and that
\be\label{Phint2}
\int_0^T\!\!\!\di s\int_0^T\!\!\!\di t\;E\big[|\Phi_i(X^1_s,X^2_t)|\big]<\infty
\qquad(T\geq 0,\ i=1,2).
\ee
Then
\be\label{apdual}
E[\Psi(X^1_T,X^2_0)]-E[\Psi(X^1_0,X^2_T)]
=\int_0^T\!\!\!\di t\;E[R(X^1_t,X^2_{T-t})]
\ee
holds for a.e.\ $T$ with respect to Lebesgue measure, where
$R(x_1,x_2):=\Phi_1(x_1,x_2)-\Phi_2(x_1,x_2) \quad(x_1\in E_1,\ x_2\in E_2)$.
Moreover, the left-hand side of (\ref{apdual}) is continuous in $T$.
\et
{\bf Proof} Although a bit of care is needed to see that all integrals are
well-defined, the proof of \cite[Theorem~7]{AS05} is correct up to the last
displayed formula (\cite[formula~(2.10)]{AS05}), which says that
\be\label{apdual2}
\int_0^S\!\di T\big(E[\Psi(X^1_T,X^2_0)]-E[\Psi(X^1_0,X^2_T)]\big)
=\int_0^S\!\di T\int_0^T\!\!\!\di t\;E[R(X^1_t,X^2_{T-t})]\qquad(S>0).
\ee
Note that by our assumption (\ref{Phint2}),
\be\label{findo}
\int_0^S\!\di T\int_0^T\!\!\!\di t\;E\big[|R(X^1_t,X^2_{T-t})|\big]
\leq\sum_{i=1}^2\int_0^S\!\!\!\di s\int_0^S\!\!\!\di t\;
E\big[|\Phi_i(X^1_s,X^2_t)|\big]<\infty
\ee
$(S>0)$, which shows that the right-hand side of (\ref{apdual2}) is
well-defined for all $S>0$. (The left-hand side of (\ref{apdual2}) is
obviously well-defined by our assumption that $\Psi\in B(E_1\times E_2)$.)
Formula (\ref{findo}) also shows that
\be
\int_0^T\!\!\!\di t\;E\big[|R(X^1_t,X^2_{T-t})|\big]<\infty
\qquad\mbox{for a.e.\ }T,
\ee
hence setting
\be
f(T):=\int_0^T\!\!\!\di t\;E[R(X^1_t,X^2_{T-t})]
\ee
yields an a.e.\ (w.r.t.\ Lebesgue measure) well-defined function $f$
satisfying $\int_0^S|f(T)|\di T<\infty$ for each $S>0$. Denoting the
left-hand side of (\ref{apdual}) by $g(T)$, formula (\ref{apdual2}) tells us
that
\be
\int_0^S\!g(T)\,\di T=\int_0^S\!f(T)\,\di T\qquad(S>0),
\ee
which implies that $g(T)=f(T)$ for a.e.\ $T$.

To finish the proof, we need to show that $T\mapsto g(T)$ is continuous. By
our assumption that $X^1$ solves the martingale problem for $A_1$,
\be
E[\Psi(X^1_T,x_2)]=\int_0^T\!\di t\,E[\Phi_1(X^1_t,x_2)]
\qquad(T\geq 0,\ x_2\in E_2).
\ee
Being an integral (which is well-defined and finite by (\ref{MP})),
the right-hand side of this equation is continuous in $T$
for each $x_2\in E_2$, hence the same is true for the left-hand side. Now if
$0\leq T_n\to T$, then by bounded pointwise convergence (using the fact that
$\Psi$ is bounded),
\be\ba{l}
\dis E[\Psi(X^1_{T_n},X^2_0)]
=\int P[X^2_0\in\di x_2]\,E[\Psi(X^1_{T_n},x_2)]\\[5pt]
\dis\qquad\asto{n}\int P[X^2_0\in\di x_2]\,E[\Psi(X^1_T,x_2)]
=E[\Psi(X^1_T,X^2_0)].
\ec
In the same way we see that $T\mapsto E[\Psi(X^1_0,X^2_T)]$ is continuous.\qed

\noi
Theorem~7 in \cite{AS05} is applied at two places in that article: in proof of
Theorem 1, pages 401--403, and in proof of Proposition 23, pages
404--405. Luckily, in both instances, all that is actually needed
is the following corollary, which still holds.
\bcor{\bf(Everywhere equality)}
Under the assumptions of Theorem~\ref{errorth}, if
\be
A_1\Psi(\,\cdot\,,x_2)(x_1)\geq A_2\Psi(x_1,\,\cdot\,)(x_2)
\qquad(x_1\in E_1,\ x_2\in E_2),
\ee
then
\be
E[\Psi(X^1_T,X^2_0)]\geq E[\Psi(X^1_0,X^2_T)]\qquad(T\in\half).
\ee
The same statement holds with both inequality signs reversed.
\ecor
{\bf Proof} Set $g(T):=E[\Psi(X^1_T,X^2_0)]-E[\Psi(X^1_0,X^2_T)]$. Then
Theorem~\ref{errorth} shows that $g$ is a continuous function satisfying
$g\geq 0$ a.e., hence $g(T)\geq 0$ for every $T\geq 0$.\qed

\noi
For completeness, we show by example that in general, the a.e.\ equality in
(\ref{apdual}) may fail to be an everywhere equality.
\begin{counterexample}
There exists metric spaces $E_i$, linear operators $A_i$ and processes $X_i$
$(i=1,2)$ together with a function $\Psi:E_1\times E_2\to\R$ satisfying the
assumptions of Theorem~\ref{errorth} such that (\ref{apdual}) does not hold
for $T=1$.
\end{counterexample}
{\bf Proof} We take $E_1=E_2=(0,\infty)$. For $r>0$, we let
$f_r:(0,\infty)\to\R$ be the function defined by $f_r(x):=e^{-rx}$. We define
linear operators $A_1,A_2$ with domains $\Di(A_1)=\Di(A_2):=\spa\{f_r:r>0\}$
by
\bc
\dis A_1f_r(x)&:=&\dis1_{\{rx\neq e\}}x\dif{x}f_r(x)
=-1_{\{rx\neq e\}}rxe^{-rx},\\[5pt]
\dis A_2f_r(x)&:=&\dis x\dif{x}f_r(x)=-rxe^{-rx}.
\ec
For $X^1,X^2$ we choose the deterministic processes
\be
X^1_t=X^2_t:=e^t\qquad(t\geq 0)
\ee
and we define $\Psi:E_1\times E_2\to\R$ by
\be
\Psi(x_1,x_2):=e^{-x_1x_2}\qquad(x_1,x_2,\geq 0).
\ee
It is straightforward to check that $X^i$ solves the martingale problem for
$A_i$ $(i=1,2)$. Note that the factor $1_{\{rx\neq e\}}$ in the definition of
$A_1$ is at this point irrelevant since for each $r>0$ there is only one time
$t$ such that $rX^1_t=e$, hence this factor has no influence on the time
integral in (\ref{duM}).

It is easy to check that (\ref{Phint2}) holds and, in the notation of
Theorem~\ref{errorth},
\be
R(x_1,x_2)=1_{\{x_1x_2=e\}}x_1x_2e^{-x_1x_2}\qquad(x_1,x_2\geq 0).
\ee
Therefore, since $X^1_tX^2_{T-t}=e^te^{T-t}=e^T$,
\be
\int_0^T\!\!\!\di t\;E[R(X^1_t,X^2_{T-t})]
=\int_0^T\!\!\!\di t\;1_{\{T=1\}}e^Te^{-e^T}
=1_{\{T=1\}}e^{1-e}\qquad(T\geq 0),
\ee
while the left-hand side of (\ref{apdual}) is in our example identically
zero.\qed

\detail{
We check that $X^i$ solves the martingale problem for $A_i$ $(i=1,2)$. Indeed,
(\ref{MP}) is easily seen to be satisfied for $f=f_r$ $(r\geq 0)$ and $X=X^i$
$(i=1,2)$. Moreover,
\be\ba{l}
\dis f_r(X^1_t)-\int_0^tA^1f_r(X^1_s)\di s
=e^{-rX^1_t}+\int_0^tr1_{\{rX^1_s\neq e\}}X^1_se^{-rX^1_s}\di s\\[5pt]
\dis\quad=e^{-re^t}+\int_0^tr1_{\{re^s\neq e\}}e^se^{-re^s}\di s\\[5pt]
\dis\quad=e^{-re^t}-\int_0^t\!\di s\,\dif{s}\big(e^{-re^s}\big)
=e^{-re^0}=e^{-r}
\ec
is a constant process, hence a martingale, and the same is true for
$f_r(X^2_t)-\int_0^tA^2f_r(X^2_s)\di s$.}

\detail{\bc
\dis\Phi_1(x_1,x_2)
&:=&\dis A_1\Psi(\,\cdot\,,x_2)(x_1)=-1_{\{x_1x_2\neq e\}}x_1x_2e^{-x_1x_2}\\[5pt]
\dis\Phi_2(x_1,x_2)
&:=&\dis A_2\Psi(x_1,\,\cdot\,)(x_2)=-x_1x_2e^{-x_1x_2},
\ec}

\begin{minipage}[t]{7cm}
{\small Siva R.~Athreya\vspace{6pt}}\\
{\small Stat.\ Math.\ Unit}\\
{\small Indian Statistical Institute}\\
{\small Bangalore Centre}\\{\small 8th Mile Mysore Road}\\
{\small Bangalore 560059, India\vspace{3pt}}\\
{\small e-mail: athreya@isibang.ac.in}\vspace{4pt}
\end{minipage}
\begin{minipage}[t]{7cm}
{\small Jan M. Swart\vspace{6pt}}\\
{\small UTIA}\\
{\small Pod vod\'arenskou v\v e\v z\' i 4}\\
{\small 18208 Praha 8}\\
{\small Czech Republic}\\
{\small e-mail: swart@utia.cas.cz} \vspace{4pt}
\end{minipage}


\begin{thebibliography}{BCGdH95}



\bibitem[AS05]{AS05}
S.R.~Athreya and J.M.~Swart.
Branching-coalescing particle systems.
{\em Prob.\ Theory Relat.\ Fields.}~131(3), 376--414, 2005.


\end{thebibliography}
\end{document}